\documentclass[11pt]{article}
\usepackage[a4paper,margin=1in]{geometry}
\usepackage{amsmath,amssymb,amsthm,mathtools}
\usepackage{enumitem}
\usepackage{booktabs}
\usepackage[hidelinks]{hyperref}
\hypersetup{colorlinks=true,linkcolor=blue,citecolor=blue,urlcolor=blue}

\newtheorem{theorem}{Theorem}[section]
\newtheorem{proposition}[theorem]{Proposition}
\newtheorem{lemma}[theorem]{Lemma}

\theoremstyle{remark}
\newtheorem{remark}[theorem]{Remark}

\newcommand{\PE}{P_E}
\newcommand{\PSL}{\mathrm{PSL}}
\newcommand{\SL}{\mathrm{SL}}
\newcommand{\PSU}{\mathrm{PSU}}
\newcommand{\SU}{\mathrm{SU}}
\newcommand{\PSp}{\mathrm{PSp}}
\newcommand{\Sp}{\mathrm{Sp}}
\newcommand{\Sz}{\mathrm{Sz}}
\newcommand{\Alt}{\mathrm{Alt}}
\newcommand{\Sym}{\mathrm{Sym}}
\newcommand{\F}{\mathbb F}
\newcommand{\sP}{\sim_P}
\newcommand{\sE}{\sim_E}
\newcommand{\nnormal}{\ntrianglelefteq}
\newcommand{\C}{C}

\title{Chordality of Enhanced Power Graphs of Finite Simple Groups: A Resolution of a Conjecture}
\author{
	Pallabi Manna
}

\date{
	\small
	Department of Mathematics, Indian Institute of Technology Hyderabad,\\ Kandi, Sangareddy-502284, Telangana, India\\
	\texttt{mannapallabimath001@gmail.com}
	\\[2mm]
	\today
}
\begin{document}
	\maketitle
	
	\begin{abstract}
		The \emph{enhanced power graph} $\PE(G)$ of a finite group $G$ has vertex set $G$, two elements being adjacent if and only if they generate a cyclic subgroup. We determine the finite simple groups $G$ for which $\PE(G)$ is chordal: they are the cyclic groups of prime order, $\Alt(7)$, $\PSL_2(q)$ ($q\ge4$), $\Sz(q)$ ($q=2^{2n+1}\ge8$) and $\PSL_3(4)$. This answers an open problem of Bubboloni, Fumagalli and Praeger \cite{BFP}. We also prove that for every finite group $G$, if the power graph $P(G)$ is a cograph (resp.\ chordal) then so is $\PE(G)$; the cograph case is a conjecture of Bubboloni, Fumagalli and Praeger \cite{BFP}.

		
	\end{abstract}
	
	\section{Introduction}\label{sec:intro}
	
	For a finite group $G$, the \emph{power graph} $P(G)$ joins $x,y$ when one lies in the cyclic group generated by the other, and the \emph{enhanced power graph} $\PE(G)$ joins $x,y$ when $\langle x,y\rangle$ is cyclic; thus $P(G)\subseteq\PE(G)$, with equality if and only if every element has prime-power order (EPPO groups) \cite{AACNS}. A graph is a \emph{cograph} if it has no induced path $P_4$, and \emph{chordal} if it has no induced cycle of length $\ge4$.
	
	\medskip
	Forbidden subgraph characterizations provide a powerful structural framework for investigating graphs defined on groups, revealing how the algebraic properties of a group govern the combinatorial structure of its associated graph. Forbidden induced subgraphs in $P(G)$ have been studied by Manna--Cameron--Mehatari \cite{MCM}, Cameron--Manna--Mehatari \cite{CMM} (simple groups with cograph power graph), Manna--Mehatari \cite{MannaMehatari2022} and Brachter--Kaja \cite{BK} (simple groups with chordal power graph), Manna--Mandal--Lucchini \cite{MML2024} (claw-free power graph); the nilpotent case for $\PE(G)$ is due to Ma--Zahirovi\'c--Lv--She \cite{MZLS}. Recently Bubboloni, Fumagalli and Praeger \cite{BFP} proved that if $\PE(G)$ is a cograph then it is chordal, classified the nonabelian simple groups with cograph $\PE(G)$ (namely $\PSL_2(q)$, $q\ge4$, $\PSL_3(4)$ and $\Sz(q)$), showed that $\PE(\Alt(7))$ is chordal but not a cograph, and posed the problem of determining the simple groups with chordal $\PE(G)$. Moreover, in \cite[\S7]{vDJMM} (see also \cite{Rundstrom}) the Kourovka problem 21.24 ($P(G)$ cograph $\Rightarrow$ $P(G)$ chordal) is solved, and the conjecture of \cite{BFP} that $P(G)$ cograph implies $\PE(G)$ cograph is recorded as open.
	
	This paper makes two contributions.
	
	\begin{theorem}[Classification]\label{thm:main}
		Let $G$ be a finite simple group. Then $\PE(G)$ is chordal if and only if $G$ is isomorphic to one of: a cyclic group of prime order; $\Alt(7)$; $\PSL_2(q)$ with $q\ge4$ (this includes $\Alt(5),\Alt(6),\PSL_3(2)$); $\Sz(q)$ with $q=2^{2n+1}\ge8$; $\PSL_3(4)$.
	\end{theorem}
	
	\begin{theorem}[Lifting]\label{thm:lift}
		Let $G$ be a finite group. If $P(G)$ is a cograph then $\PE(G)$ is a cograph; if $P(G)$ is chordal then $\PE(G)$ is chordal.
	\end{theorem}
	
	\noindent Together with \cite{BFP} and \cite{vDJMM,Rundstrom} we obtain, for every finite group $G$,
	\[
	P(G)\ \text{cograph}\Longrightarrow \PE(G) \text{cograph}\Longrightarrow\ \PE(G) \text{chordal},
	P(G) \text{cograph}\Longrightarrow\ P(G) \text{chordal}\Longrightarrow \PE(G)\text{chordal}.
	\]
	
	\paragraph{The principal contribution:} By \cite[Theorem~B]{BFP}, the groups $\PSL_2(q)$, $\PSL_3(4)$, $\Sz(q)$ have cograph, hence chordal, $\PE(G)$, and $\Alt(7)$ is chordal; so the \emph{sufficiency} half of Theorem~\ref{thm:main} is already contained in \cite{BFP} (we include direct proofs for completeness; $\Alt(5),\Alt(6),\Alt(7),\PSL_3(4)$ also follow from \cite[Theorem~1.1]{BK} and Theorem~\ref{thm:lift}). The new content is the \emph{necessity}: every other simple group has non-chordal $\PE(G)$, which answers the chordal part of Open Problem~1 of \cite{BFP}; in particular $\Alt(7)$ is the only simple group with chordal but non-cograph $\PE(G)$. Theorem~\ref{thm:lift} proves, for all finite groups, the cograph conjecture of \cite{BFP}, and its chordal version. The converses fail: $\PE(C_{30})$ is complete, whereas $P(C_{30})$ is not chordal (the elements of orders $6,2,10,5,15,3$ induce a $C_6$, since adjacency in $P(C_{30})$ is divisibility of orders); and $P(\mathbb Z_{12})$ is not a cograph (elements $3,6,2,4$ of orders $4,2,6,3$ induce a $P_4$) while $\PE(\mathbb Z_{12})$ is complete. Within simple groups, $\PSL_2(59)$ contains a cyclic subgroup of order $30$, so $P(\PSL_2(59))$ is not chordal, while $\PE(\PSL_2(59))$ is chordal (Prop.~\ref{prop:psl2}).
	
	The negative half rests on one mechanism (Lemma~\ref{lem:E}): commuting elements $x,y$ of coprime prime-power orders with neither $\langle x\rangle$ normal in $\C_G(y)$ nor $\langle y\rangle$ normal in $\C_G(x)$ yield an induced even cycle; the positive half on a perfect elimination ordering by vertices lying in a unique maximal cyclic subgroup.\\
	\textbf{Organization:} Section~\ref{sec:tools} collects fundamental concepts and results, Section~\ref{sec:lift} proves Theorem~\ref{thm:lift}, Sections~\ref{sec:pos} and~\ref{sec:neg} treat chordal and non-chordal families, Section~\ref{sec:proof} assembles the proof of Theorem~\ref{thm:main}, and Section~\ref{sec:comp} records computer checks using GAP~4.12. Section \ref{sec:facts} summarizes some well-known results and standard facts.
	
	\section{Preliminary Results}\label{sec:tools}
	
	Throughout, we use the symbol $x^c$ to mean $x^c=cxc^{-1}$. We write $x\sP y$ if $x\ne y$ and ($x\in\langle y\rangle$ or $y\in\langle x\rangle$), and $x\sE y$ if $x\neq y$ and $\langle x,y\rangle$ is cyclic. If $x\in\langle y\rangle$ then $\langle x,y\rangle=\langle y\rangle$, so $\sP\,\subseteq\,\sE$, and every non-edge of $\PE(G)$ is a non-edge of $P(G)$. We write $A\nnormal B$ for ``$A$ is not normal in $B$''.
	
	\begin{lemma}[Subgroups and quotients]\label{lem:sub}
		(a) If $H\le G$ then $\PE(H)$ is the induced subgraph of $\PE(G)$ on $H$; so if $H$ is not chordal, neither is $G$.\\
		(b) If $N\trianglelefteq G$, $H\le G$ and $H\cap N=1$, then $H\cong HN/N\le G/N$; so if $H$ is not chordal, neither is $G/N$.
	\end{lemma}
	\begin{proof}
		Adjacency depends only on $\langle x,y\rangle$, giving (a); (b) is the isomorphism theorem and (a).
	\end{proof}
	
	\begin{lemma}[Universal vertex]\label{lem:univ}
		The identity is adjacent to every vertex, so it lies on no induced cycle of length $\ge4$; $\PE(G)$ is chordal if and only if $\PE(G)-\{1\}$ is.
	\end{lemma}
	
	\begin{lemma}[EPPO groups \cite{AACNS}]\label{lem:eppo}
		If every element of $G$ has prime-power order, then $\PE(G)=P(G)$ and it is chordal.
	\end{lemma}
	\begin{proof}
		If $\langle x,y\rangle$ is cyclic of prime-power order its subgroups form a chain, so $x\sP y$ or $x=y$. Let $C$ be an induced cycle of length $\ge4$ and $x\in C$ of maximal order; its neighbours on $C$ are comparable with $x$ and have order $\le|x|$, so lie in $\langle x\rangle$; subgroups of a cyclic $p$-group form a chain, so they are adjacent, a contradiction.
	\end{proof}
	
	\begin{lemma}[Simplicial elimination \cite{FG}]\label{lem:simp}
		A vertex is \emph{simplicial} if its closed neighbourhood is a clique. If the vertices can be ordered $v_1,\dots,v_n$ so that each $v_i$ is simplicial in the subgraph induced on $\{v_i,\dots,v_n\}$, the graph is chordal. In $\PE(G)$ the closed neighbourhood of $x$ is the union of the cyclic subgroups containing $x$; if $x$ lies in a unique maximal cyclic subgroup $M$ then $N[x]=M$ and $x$ is simplicial.
	\end{lemma}
	\begin{proof}
		Let $C$ be an induced cycle of length $\ge4$ and $v_i$ its first vertex in the order. Its two neighbours on $C$ lie in $\{v_i,\dots,v_n\}$ and are non-adjacent, so $v_i$ is not simplicial there.
	\end{proof}
	
	\begin{lemma}[Coprime commuting pairs]\label{lem:E}
		Let $x,y\in G$ commute and have coprime prime-power orders $r^a,s^b$. If $\langle x\rangle\nnormal \C_G(y)$ and $\langle y\rangle\nnormal \C_G(x)$, then $\PE(G)$ has an induced cycle of even length $\ge4$.
	\end{lemma}
	\begin{proof}
		Let $\mathcal X,\mathcal Y$ be the sets of $G$-conjugates of $\langle x\rangle,\langle y\rangle$ and choose one generator for each member; these are the vertices of $B$. Two distinct members of $\mathcal X$ are non-adjacent because a cyclic group has at most one subgroup of order $r^a$; likewise for $\mathcal Y$. A member of $\mathcal X$ and one of $\mathcal Y$ are adjacent if and only if they commute (commuting elements of coprime orders generate a cyclic group). So $B$ is a bipartite induced subgraph of $\PE(G)$. The vertex $\langle x\rangle$ is adjacent to every $\langle y\rangle^c$ with $c\in \C_G(x)$, and there are $|\C_G(x):N_{\C_G(x)}(\langle y\rangle)|\ge2$ of these; by conjugation every vertex of $\mathcal X$ has degree $\ge2$, and symmetrically for $\mathcal Y$. A finite graph of minimum degree $\ge2$ has a cycle, a shortest cycle is induced, and in $B$ it has even length $\ge4$.
	\end{proof}
	
	\begin{lemma}[Quotient version]\label{lem:Eq}
		Let $Z\le Z(G)$ and $\bar G=G/Z$. Let $x,y\in G$ commute with $\bar x,\bar y$ of coprime prime-power orders. If there are $c\in \C_G(y)$ and $c'\in \C_G(x)$ with $x^c\notin\langle x\rangle Z$ and $y^{c'}\notin\langle y\rangle Z$, then $\PE(\bar G)$ has an induced cycle of length $\ge4$.
	\end{lemma}
	\begin{proof}
		$\bar c\in \C_{\bar G}(\bar y)$ and $\bar x^{\bar c}\notin\langle\bar x\rangle$, so $\langle\bar x\rangle\nnormal \C_{\bar G}(\bar y)$; symmetrically for $\bar y$. Apply Lemma~\ref{lem:E}.
	\end{proof}
	
	\begin{lemma}[Numerical criterion]\label{lem:num}
		Let $x,y$ commute and have prime orders $r\neq s$, and put $g=\gcd(|\C_G(x)|,|\C_G(y)|)$. If $|\C_G(x)|/g>s-1$ and $|\C_G(y)|/g>r-1$, then $\PE(G)$ is not chordal.
	\end{lemma}
	\begin{proof}
		If $\langle y\rangle\trianglelefteq \C_G(x)$ then $\C_G(x)/(\C_G(x)\cap \C_G(y))$ embeds in $\mathrm{Aut}(C_s)$, so $|\C_G(x):\C_G(x)\cap \C_G(y)|\le s-1$. But $|\C_G(x)\cap \C_G(y)|$ divides $g$, so this index is $\ge|\C_G(x)|/g$, a contradiction. The same holds for $\langle x\rangle$. Apply Lemma~\ref{lem:E}.
	\end{proof}
	
	\begin{lemma}[Direct products]\label{lem:prod}
		Suppose $A\times B\le G$, where $A$ has two distinct subgroups of prime order $r$ and $B$ has two distinct subgroups of prime order $s\ne r$. Then $\PE(G)$ is not chordal. This applies when $A,B$ are non-abelian simple, with $r\mid|A|$, $s\mid|B|$, $r\ne s$ (a subgroup of prime order in a simple group is not normal, hence has $\ge2$ conjugates), and to $\Alt(4)\times\Alt(5)$ and $\PSL_3(2)\times S_4$.
	\end{lemma}
	\begin{proof}
		Take $a,a'\in A$ of order $r$ and $b,b'\in B$ of order $s$ with $\langle a\rangle\ne\langle a'\rangle$, $\langle b\rangle\ne\langle b'\rangle$. Then $a\not\sE a'$ and $b\not\sE b'$, while every cross pair commutes and has coprime orders, hence is adjacent. So $(a,b,a',b')$ is an induced $C_4$.
	\end{proof}
	
	\section{The lifting theorem}\label{sec:lift}
	
	\begin{proof}[Proof of Theorem~\ref{thm:lift}]
		We prove the contrapositive: if $\PE(G)$ contains an induced cycle of length $\ge4$ (resp.\ an induced $P_4$), so does $P(G)$.
		
		Let $\Gamma$ be an induced subgraph of $\PE(G)$ that is a cycle $v_1,\dots,v_k$ ($k\ge4$, indices mod $k$) or a path $v_1,\dots,v_k$ ($k\ge4$), with edges $e_i=v_iv_{i+1}$. Being induced and of length $\ge4$, $\Gamma$ has no triangle.
		
		\smallskip\noindent\emph{\textbf{Step 1}.} Let $I$ be the set of $i$ such that $e_i$ is an edge of $\Gamma$ but $v_i\not\sP v_{i+1}$. For $i\in I$ let $z_i$ generate the cyclic group $\langle v_i,v_{i+1}\rangle$; then $v_i,v_{i+1}\in\langle z_i\rangle$.
		
		\smallskip\noindent\emph{\textbf{Step 2} (key observation (K)).} Let $i\in I$ and let $w$ be a vertex of $\Gamma$. If $w\in\langle z_i\rangle$ or $z_i\in\langle w\rangle$, then $w\in\{v_i,v_{i+1}\}$.
		
		Indeed, if $w\in\langle z_i\rangle$ then $\langle w,v_i\rangle$ and $\langle w,v_{i+1}\rangle$ lie in the cyclic group $\langle z_i\rangle$; if $z_i\in\langle w\rangle$ then $v_i,v_{i+1}\in\langle w\rangle$ and both groups equal $\langle w\rangle$. In either case $w$ equals or is $\sE$-adjacent to each of $v_i,v_{i+1}$. If $w\notin\{v_i,v_{i+1}\}$ this is a triangle in $\Gamma$, which is impossible.
		
		\smallskip\noindent\emph{\textbf{Step 3} (consequences).}
		\begin{enumerate}[label=(\alph*)]
			\item $z_i\notin V(\Gamma)$: if $z_i=v_j$ then $v_j\in\langle z_i\rangle$, so by (K) $v_j\in\{v_i,v_{i+1}\}$; but $z_i\in\{v_i,v_{i+1}\}$ would give $v_i\sP v_{i+1}$, contradicting $i\in I$.
			\item For $i\ne j$ in $I$, $z_j\notin\langle z_i\rangle$ (in particular $z_i\neq z_j$ and $z_i\not\sP z_j$): if $z_j\in\langle z_i\rangle$ then $v_j,v_{j+1}\in\langle z_j\rangle\le\langle z_i\rangle$, so by (K) $\{v_j,v_{j+1}\}=\{v_i,v_{i+1}\}$ and $j=i$.
			\item For $i\in I$ and a vertex $v_j$ of $\Gamma$: $z_i\sP v_j$ if and only if $j\in\{i,i+1\}$. The ``if'' holds as $v_j\in\langle z_i\rangle$; the ``only if'' is (K).
		\end{enumerate}
		
		\smallskip\noindent\emph{\textbf{Step 4.}} Let $\Gamma'$ arise from $\Gamma$ by replacing, for every $i\in I$, the edge $v_iv_{i+1}$ by the path $v_i,z_i,v_{i+1}$. By (a),(b) its vertices are distinct, and it is a cycle (resp.\ path) of length $\ge k\ge4$. We claim $\Gamma'$ is an induced subgraph of $P(G)$.
		\begin{itemize}
			\item Its edges are $\sP$-edges: $v_iz_i$, $z_iv_{i+1}$ by (c); $e_i$ for $i\notin I$ by definition of $I$.
			\item A pair $v_jv_l$ that is not an edge of $\Gamma'$ is non-adjacent in $\Gamma$ (hence in $\PE(G)$, hence in $P(G)$), or equals $v_iv_{i+1}$ with $i\in I$ (non-adjacent in $P(G)$ by definition of $I$).
			\item A pair $z_iv_j$ that is not an edge has $j\notin\{i,i+1\}$ and is non-adjacent by (c).
			\item $z_iz_j$ ($i\ne j$) is non-adjacent by (b).
		\end{itemize}
		
		\smallskip\noindent\emph{\textbf{Step 5.}} If $\Gamma$ is an induced cycle, $\Gamma'$ is an induced cycle of length $\ge4$ in $P(G)$. If $\Gamma$ is an induced $P_4$, $\Gamma'$ is an induced path on $\ge4$ vertices of $P(G)$ and any four consecutive vertices induce a $P_4$.
	\end{proof}
	
	\section{Simple groups with chordal $\PE(G)$}\label{sec:pos}
	
	The following facts are used; sources are listed in Section~\ref{sec:facts}.
	
	\begin{proposition}[$\PSL_2(q)$, $q\ge4$]\label{prop:psl2}
		$\PE(\PSL_2(q))$ is chordal, and in fact $\PE(\PSL_2(q))-\{1\}$ is a disjoint union of cliques; in particular $\PE(\PSL_2(q))$ is a cograph (cf.\ \cite{BFP}).
	\end{proposition}
	\begin{proof}
		Facts D (Section~\ref{sec:facts}). Let $A$ be a Sylow $p$-subgroup ($q=p^f$); $\C(a)=A$ for $a\in A^\#$.
		
		\emph{$q$ even.} All involutions lie in Sylow $2$-subgroups, and the tori have odd orders $q\pm1$. Every $u\ne1$ is a $2$-element, in $A$, or lies in a torus $T$ with $\C(u)=T$.
		
		\emph{$q$ odd.} For an involution $z$, $\C(z)$ is dihedral of order $q-\varepsilon$ ($q\equiv\varepsilon\bmod 4$), $z$ is its central involution and the rotation subgroup $R$ is cyclic of even order $m=(q-\varepsilon)/2$.
		
		We show every $u\ne1$ lies in a unique maximal cyclic subgroup $M(u)$:
		\begin{itemize}
			\item $u\in A$: a cyclic subgroup containing $u$ lies in $\C(u)=A$, so has order $p$ and $M(u)=\langle u\rangle$.
			\item $|u|>2$ and $p\nmid|u|$: a cyclic subgroup containing $u$ lies in $\C(u)=T$, so $M(u)=T$.
			\item $u=z$ an involution, $q$ odd: a cyclic subgroup containing $z$ lies in $\C(z)$. If $m\ge4$ a reflection is not central, so it lies in $R$ and $M(z)=R$. If $m=2$ ($q=5$) then $\C(z)\cong V_4$ and $M(z)=\langle z\rangle$.
		\end{itemize}
		If $\langle u,h\rangle$ is cyclic it lies in a maximal cyclic subgroup, which then equals $M(u)$ and $M(h)$; so $u\sE h\iff M(u)=M(h)$, and $\PE-\{1\}$ is a disjoint union of cliques. Apply Lemma~\ref{lem:univ}.
	\end{proof}
	
	\begin{proposition}[$\Sz(q)$]\label{prop:sz}
		Let $q=2^{2n+1}$ with $n\ge1$. Then $\PE(\Sz(q))$ is chordal.
	\end{proposition}
	\begin{proof}
		Facts S. Let $P$ be a Sylow $2$-subgroup, $Z=Z(P)=\Omega_1(P)$.
		
		\emph{Odd order.} If $|u|$ is odd, neighbours of $u$ lie in $\C(u)=T$, a cyclic torus, so $N[u]=T$ is a clique.
		
		\emph{Involution $z\in Z$} (note $\C(z)=P$ for every involution $z\in Z$, and all involutions lie in $Z$). A neighbour $h$ lies in $P$. If $|h|=2$ and $h\ne z$ then $\langle z,h\rangle\cong V_4$ is not cyclic. If $|h|=4$ then $h^2$ is the unique involution of $\langle z,h\rangle$, so $h^2=z$.
		
		\emph{Element $w$ of order $4$.} A neighbour lies in $\C(w)=\langle w\rangle Z$, so is $h=w^iz'$ with $z'\in Z$ and $\langle w,h\rangle=\langle w,z'\rangle$; if $z'\notin\langle w^2\rangle$ this contains $V_4$. So $N[w]=\langle w\rangle$.
		
		Hence $\PE(\Sz(q))-\{1\}$ is the disjoint union of the cliques given by the tori (minus $1$) and the components $\{z\}\cup\{h:h^2=z\}$, each a union of triangles $\{z,w,w^{-1}\}$ sharing the cut vertex $z$. Every block is a clique, so there is no induced cycle of length $\ge4$; apply Lemma~\ref{lem:univ}. (Note that $\PE(\Sz(q))$ itself contains diamonds $\{1,z,w,w'\}$ with $w^2=w'^2=z$, $\langle w\rangle\ne\langle w'\rangle$; the block statement refers to $\PE-\{1\}$.)
	\end{proof}
	
	\begin{proposition}[$\PSL_3(4)$]\label{prop:psl34}
		$\PE(\PSL_3(4))$ is chordal.
	\end{proposition}
	\begin{proof}
		The element orders of $\PSL_3(4)$ are $1,2,3,4,5,7$ (ATLAS), so it is EPPO; apply Lemma~\ref{lem:eppo}.
	\end{proof}
	
	\begin{proposition}[$\Alt(7)$]\label{prop:a7}
		$\PE(\Alt(7))$ is chordal.
	\end{proposition}
	\begin{proof}
		The element orders are $1,\dots,7$ with $4,5,6,7$ occurring. A cyclic overgroup of an element of order $4,5,6$ or $7$ has order a multiple of it and at most $7$, so equal to it; thus such elements lie in a unique maximal cyclic subgroup, are simplicial (Lemma~\ref{lem:simp}), and may be removed first.
		
		All involutions of $\Alt(7)$ are conjugate to $z=(12)(34)$, and $\C_{\Alt(7)}(z)=\tfrac12|\C_{\Sym(7)}(z)|=\tfrac12|D_8\times\Sym\{5,6,7\}|=24$ has Sylow $3$-subgroup $\langle(567)\rangle$ (a subgroup of order $3$ of $D_8\times\Sym\{5,6,7\}$ lies in $\Sym\{5,6,7\}$, where it is unique). In the remaining graph $N[z]=\{1,z,u,u^2\}$ with $u=(567)$, which is a clique ($\langle zu\rangle$ is cyclic of order $6$); so $z$ is simplicial and all involutions may be removed. What remains is $\{1\}$ together with the subgroups of order $3$, triangles through $1$; these eliminate in turn. By Lemma~\ref{lem:simp}, $\PE(\Alt(7))$ is chordal.
	\end{proof}
	
	\section{Simple groups with non-chordal $\PE(G)$}\label{sec:neg}
	
	\subsection{Alternating groups}
	
	\begin{proposition}\label{prop:alt}
		$\PE(\Alt(n))$ is non-chordal for $n\ge8$. Also $\PE(\Sym(n))$ is non-chordal for $n\ge6$ (only $\Sym(6)\cong\Sp_4(2)$ is used below, as an auxiliary group).
	\end{proposition}
	\begin{proof}
		$\Sym(6)$: take $x=(123)$, $y=(45)$. Then $\C(x)=\langle x\rangle\times\Sym\{4,5,6\}$ contains $(46)$, which moves $\langle y\rangle$, and $\C(y)=\langle y\rangle\times\Sym\{1,2,3,6\}$ contains $(16)$ with $x^{(16)}=(623)\notin\langle x\rangle$. Lemma~\ref{lem:E}.
		
		$\Alt(8)$: take $x=(123)$, $y=(45)(67)$. Then $\C_{\Alt(8)}(x)=\langle x\rangle\times\Alt\{4,\dots,8\}$ contains $y\in\Alt(5)$ and $\langle y\rangle$ is not normal in the simple group $\Alt(5)$. Also $c=(18)(45)\in\Alt(8)$ commutes with $y$ and $x^c=(823)\notin\langle x\rangle$. Lemma~\ref{lem:E}.
		
		Larger $n$: Lemma~\ref{lem:sub}(a).
	\end{proof}
	
	\subsection{Linear groups}
	
	Let $E_{ij}$ denote matrix units.
	
	\begin{lemma}[Central quotients of $\SL_3(q)$]\label{lem:sl3}
		Let $q\ne2,4$ and $Z_0\le Z(\SL_3(q))$. Then $\PE(\SL_3(q)/Z_0)$ is non-chordal. (In particular $\SL_3(q)$ and $\PSL_3(q)$ are non-chordal.)
	\end{lemma}
	\begin{proof}
		Put $y=I+E_{12}$, $c=I+E_{13}$, $c'=I+E_{21}$ and choose $x$ as follows.
		\begin{itemize}
			\item $q$ odd: $x=\mathrm{diag}(-1,-1,1)$, of order $2$; $x\notin Z_0$ since $|Z(\SL_3(q))|$ divides $3$, so $\bar x$ has order $2$; $|\bar y|=p$ is odd.
			\item $q=2^f\ge8$: $q-1\equiv7\pmod 8$ while $3^k\equiv1,3\pmod 8$, so $q-1$ has a prime divisor $r\ne3$. Take $a\in\F_q^\times$ of order $r$ and $x=\mathrm{diag}(a,a,a^{-2})$. Then $a^3\ne1$, so $x$ is not scalar, $\bar x$ has order $r$, and $|\bar y|=2$.
		\end{itemize}
		Then $x$ commutes with $y$ and $c'$ (its first two diagonal entries are equal), and $c\in \C(y)$ because $E_{12}E_{13}=E_{13}E_{12}=0$. Further $cxc^{-1}=x+(x_{33}-x_{11})E_{13}$ with $x_{33}\ne x_{11}$ is not diagonal, so $x^c\notin\langle x\rangle Z_0$, a group of diagonal matrices. And $y^{c'}=I-E_{11}+E_{22}+E_{12}-E_{21}$ has non-zero $(2,1)$ entry, whereas every element of $\langle y\rangle Z_0$ is a scalar multiple of $I+kE_{12}$. Lemma~\ref{lem:Eq} applies.
	\end{proof}
	
	\begin{proposition}\label{prop:psln}
		$\PE(\PSL_n(q))$ is non-chordal for $n=3$, $q\ne2,4$, and for all $n\ge4$ and all $q$.
	\end{proposition}
	\begin{proof}
		$n=3$ is Lemma~\ref{lem:sl3}. Let $n\ge4$ and $Z=\{\lambda I:\lambda^n=1\}$.
		\begin{itemize}
			\item $q\ne2,4$: the block $H=\SL_3(q)\oplus I_{n-3}$ meets $Z$ trivially (a scalar in $H$ is $1$ on the identity block); Lemmas~\ref{lem:sl3} (with $Z_0=1$) and~\ref{lem:sub}(b).
			\item $q=4$: $H=(\SL_2(4)\times\SL_2(4))\oplus I_{n-4}\cong\Alt(5)\times\Alt(5)$ meets $Z$ trivially (a scalar $\lambda$ in a block $\SL_2(4)$ has $\lambda^2=1$, so $\lambda=1$); Lemmas~\ref{lem:prod},~\ref{lem:sub}(b).
			\item $q=2$: $\PSL_n(2)=\SL_n(2)\supseteq\SL_4(2)\cong\Alt(8)$; Prop.~\ref{prop:alt}.\qedhere
		\end{itemize}
	\end{proof}
	
	\subsection{Symplectic groups}
	
	Let $e_1,f_1,e_2,f_2$ be a symplectic basis with $B(e_i,f_i)=1=-B(f_i,e_i)$, and $T_u:v\mapsto v+B(v,u)u$ a symplectic transvection; then $gT_ug^{-1}=T_{gu}$ and $T_{-u}=T_u$.
	
	\begin{proposition}\label{prop:sp4}
		For $q$ odd, $\Sp_4(q)$ and $\PSp_4(q)$ are non-chordal; for $q$ even, $q\ge4$, $\Sp_4(q)=\PSp_4(q)$ is non-chordal; $\Sp_4(2)\cong\Sym(6)$ is non-chordal. (For $q=2$ the simple group is $\Sp_4(2)'\cong\Alt(6)\cong\PSL_2(9)$, which \emph{is} chordal.)
	\end{proposition}
	\begin{proof}
		$q=2$: Prop.~\ref{prop:alt}.
		
		\emph{$q$ odd.} Let $x=-1$ on $W_1=\langle e_1,f_1\rangle$ and $+1$ on $W_2=\langle e_2,f_2\rangle$; $y=T_{e_1}$, $c=T_{e_1+e_2}$, $c'=T_{f_1}$, $Z=\{\pm1\}$. Since $xT_ux^{-1}=T_{xu}$ and $T_{-u}=T_u$ we get $[x,y]=[x,c']=1$; $ce_1=e_1$ as $B(e_1,e_1+e_2)=0$, so $c\in \C(y)$. Further $cf_1=f_1+B(f_1,e_1+e_2)(e_1+e_2)=f_1-(e_1+e_2)\notin W_1\cup W_2$, so the $(-1)$-eigenspace $c(W_1)$ of $x^c$ is neither $W_1,W_2$ nor $0,V$; hence $x^c\notin\{\pm1,\pm x\}$. Finally $c'e_1=e_1+B(e_1,f_1)f_1=e_1+f_1$, so $y^{c'}=T_{e_1+f_1}$ is a non-trivial unipotent element with image line $\langle e_1+f_1\rangle\ne\langle e_1\rangle$, while the unipotent elements of $\langle y\rangle Z$ are $1$ or have image line $\langle e_1\rangle$. So $y^{c'}\notin\langle y\rangle Z$. The elements $x,y$ have orders $2,p$. By Lemma~\ref{lem:Eq} with $Z$, $\PSp_4(q)$ is non-chordal. The same computation shows $x^c\notin\{1,x\}$ and $y^{c'}\notin\langle y\rangle$, so Lemma~\ref{lem:Eq} with $Z=1$ shows that $\Sp_4(q)$ itself is non-chordal.
		
		\emph{$q$ even, $q\ge4$.} $\Sp_4(q)\ge\SL_2(q)\times\SL_2(q)$ with simple factors, and the centre is trivial. Lemma~\ref{lem:prod} with $r=2$ and an odd prime $s\mid q+1$.
	\end{proof}
	
	\begin{proposition}\label{prop:spbig}
		$\PE(\PSp_{2m}(q))$ is non-chordal for all $m\ge3$ and all $q$.
	\end{proposition}
	\begin{proof}
		$H=\Sp_4(q)\oplus I_{2m-4}$ meets the centre $\{\pm1\}$ trivially (as $-1$ is non-trivial on the identity block, of dimension $2m-4\ge2$). By Prop.~\ref{prop:sp4}, $H\cong\Sp_4(q)$ is non-chordal for every $q$ (including $q=2$, where $\Sp_4(2)\cong\Sym(6)$); apply Lemma~\ref{lem:sub}(b).
	\end{proof}
	
	\subsection{Unitary groups}
	
	We use the Hermitian form with Gram matrix $J=\mathrm{antidiag}(1,1,1)$ on $\F_{q^2}^3$; bars denote $\lambda\mapsto\lambda^q$. The direct computation of $\bar g^{T}Jg=J$ gives:
	\begin{itemize}
		\item $\begin{psmallmatrix}1&\alpha&\beta\\0&1&-\bar\alpha\\0&0&1\end{psmallmatrix}\in\SU_3(q)$ if and only if $\beta+\bar\beta+\alpha\bar\alpha=0$; these form a group $U$ of order $q^3$ with centre $\{I+tE_{13}:t+\bar t=0\}$;
		\item $I+sE_{31}\in\SU_3(q)$ if and only if $s+\bar s=0$;
		\item $\mathrm{diag}(a,b,a')\in\SU_3(q)$ if and only if $b\bar b=1$, $a'=\bar a^{-1}$ and $aba'=1$.
	\end{itemize}
	
	\begin{proposition}\label{prop:su3}
		$\PE(\SU_3(q))$ and $\PE(\PSU_3(q))$ are non-chordal for $q\ge3$.
	\end{proposition}
	\begin{proof}
		Let $Z=Z(\SU_3(q))$, of order $\gcd(3,q+1)$; the case $Z=1$ gives $\SU_3(q)$. Choose $a\in\F_{q^2}^\times$ with $a^{q+1}=1$, $a^3\ne1$ and of prime-power order:
		\begin{itemize}
			\item $q$ odd: $a=-1$;
			\item $q$ even, $q+1$ has a prime divisor $r\ne3$: $a$ of order $r$;
			\item $q$ even, $q+1=3^k$: by Mih\u{a}ilescu's theorem $q=8$; take $a$ of order $9$ (then $x^3=a^3I\in Z$ and $\bar x$ has order $3$).
		\end{itemize}
		In every case $a\ne a^{-2}$ (as $a^3\ne1$), so $x=\mathrm{diag}(a,a^{-2},a)\in\SU_3(q)$ is not scalar and $\bar x\ne1$. Put $y=I+tE_{13}$ with $t\ne0$, $t+\bar t=0$; $|y|=p$ is coprime to $|\bar x|$ (for $q$ odd $|x|=2$; for $q$ even $|y|=2$ and $|\bar x|$ is odd). Clearly $[x,y]=1$.
		
		Take $c\in U$ with $\alpha\ne0$ (exists since the trace $\F_{q^2}\to\F_q$ is onto). As $E_{13}$ is central in $U$, $c\in \C(y)$. The $(1,2)$ entry of $cxc^{-1}-x$ is $\alpha(a^{-2}-a)\ne0$, so $x^c$ is not diagonal and $x^c\notin\langle x\rangle Z$.
		
		Take $c'=I+sE_{31}$ with $s\ne0$, $s+\bar s=0$; it commutes with $x$ since $x_{11}=x_{33}$, and
		\[
		y^{c'}=I+t\,(E_{13}-sE_{11}+sE_{33}-s^2E_{31})
		\]
		has $(3,1)$ entry $-ts^2\ne0$, while every element of $\langle y\rangle Z$ has $(3,1)$ entry $0$. Apply Lemma~\ref{lem:Eq}.
	\end{proof}
	
	\begin{proposition}\label{prop:sun}
		$\PE(\PSU_n(q))$ is non-chordal for all $n\ge4$ and all $q$.
	\end{proposition}
	\begin{proof}
		$q\ge3$: the block $\SU_3(q)\oplus I_{n-3}$ meets the centre (scalars) trivially; Prop.~\ref{prop:su3}, Lemma~\ref{lem:sub}(b). $q=2$, $n=4$: $Z(\SU_4(2))=1$ and $\PSU_4(2)\cong\PSp_4(3)$, non-chordal by Prop.~\ref{prop:sp4}. $q=2$, $n\ge5$: block $\SU_4(2)\oplus I_{n-4}$.
	\end{proof}
	
	\subsection{Orthogonal groups}
	
	Let $O$ be the full isometry group of the quadratic form, and $\Omega\le O$ the usual subgroup: for $n\ge5$ the derived subgroup of $O$; equivalently the kernel of the spinor norm and determinant ($q$ odd), resp.\ of the Dickson invariant ($q$ even). In all cases $O/\Omega$ is abelian (of order $\le4$), so $\Omega$ contains every perfect subgroup of $O$.
	
	\begin{proposition}\label{prop:orth}
		Let $G$ be $\Omega_{2n+1}(q)$ ($n\ge3$), $P\Omega^{+}_{2n}(q)$ ($n\ge4$) or $P\Omega^-_{2n}(q)$ ($n\ge4$). Then $\PE(G)$ is non-chordal. (Orthogonal groups of lower rank are isomorphic to $\PSL_2(q)$, $\PSL_2(q^2)$, $\PSp_4(q)$, $\PSL_4(q)$, $\PSU_4(q)$ and are covered elsewhere.)
	\end{proposition}
	\begin{proof}
		\emph{$q\ne2,4$.} Write the space as $U\oplus U^*\oplus W$ with $U$ totally singular of dimension $k\ge3$ ($k=n$, $\dim W=1$ for $\Omega_{2n+1}$; $k=n$, $W=0$ for the $+$ type; $k=n-1$, $W$ an anisotropic $2$-space, for the $-$ type). The map $g\mapsto g\oplus g^{-T}\oplus1$ embeds $\SL_3(q)\le\SL_k(q)$ (acting on the first three coordinates of $U$) into $O$; the image $H\cong\SL_3(q)$ is perfect, so $H\le\Omega$. $H$ acts trivially on a non-zero subspace, so $H$ meets the centre of $\Omega$ ($\subseteq\{\pm1\}$) trivially. By Lemma~\ref{lem:sl3} with $Z_0=1$ and Lemma~\ref{lem:sub}(b), $G$ is non-chordal. (For $q$ odd and type $2n+1$ the centre is trivial.)
		
		\emph{$q\in\{2,4\}$.} For $\Omega_{2n+1}(q)\cong\Sp_{2n}(q)$ use Prop.~\ref{prop:spbig}. For $\Omega^\pm_{2n}(q)$, the stabilizer in $O^\pm_{2n}(q)$ of a non-singular vector $v$ is $\Sp_{2n-2}(q)\times\langle r_v\rangle$ ($r_v$ the transvection with centre $v$, which fixes $v$ as $q$ is even). The factor $\Sp_{2n-2}(q)$ is perfect since $2n-2\ge6$, so it lies in $\Omega$; it is non-chordal by Prop.~\ref{prop:spbig}. The centre is trivial for $q$ even.
	\end{proof}
	
	\subsection{Exceptional groups and the Tits group}
	
	\begin{proposition}\label{prop:exc}
		$\PE(G)$ is non-chordal for $G={}^2G_2(q)$ ($q=3^{2n+1}\ge27$), $G_2(q)$ ($q\ge3$), $G_2(2)'\cong\PSU_3(3)$, ${}^3D_4(q)$, $F_4(q)$, ${}^2F_4(q)$ ($q=2^{2n+1}\ge8$), the Tits group ${}^2F_4(2)'$, $E_6(q)$, ${}^2E_6(q)$, $E_7(q)$, $E_8(q)$.
	\end{proposition}
	\begin{proof}
		\emph{${}^2G_2(q)$, $q\ge27$.} By Facts W, an involution $x$ has $|\C(x)|=q(q^2-1)$ with $\C(x)\cong C_2\times\PSL_2(q)$; $\C(x)$ contains an element $y$ of order $3$ (as $3\mid q$), and $|\C(y)|\in\{q^3,2q^2\}$. Then $g=\gcd(|\C(x)|,|\C(y)|)\in\{q,2q\}$, so $|\C(x)|/g\ge(q^2-1)/2>2$ and $|\C(y)|/g\ge q>1$. Lemma~\ref{lem:num}.
		
		\emph{Containments} (sources in Table~\ref{tab:cont}): $G_2(q)\supseteq\SU_3(q)$ for $q\ge3$, and $G_2(2)'\cong\PSU_3(3)$ (ATLAS); ${}^3D_4(q)\supseteq G_2(q)$; the Tits group contains $L_3(3){:}2\supseteq\SL_3(3)=\PSL_3(3)$ (non-chordal by Lemma~\ref{lem:sl3} as $3\ne2,4$) and ${}^2F_4(q)\supseteq{}^2F_4(2)'$ (subfield subgroup); $F_4(q)$ contains the subgroup generated by long-root subgroups of a subsystem of type $A_2$, which is $\SL_3(q)/Z_0$ for some $Z_0\le Z(\SL_3(q))$, non-chordal for $q\ne2,4$ by Lemma~\ref{lem:sl3}, whichever $Z_0$ occurs; for $q$ even $F_4(q)\supseteq B_4(q)\cong\Sp_8(q)$, non-chordal by Prop.~\ref{prop:spbig}, which covers $q=2,4$. For $G\in\{E_6(q),{}^2E_6(q),E_7(q),E_8(q)\}$ (simple versions), let $\tilde G$ be the corresponding simply connected group; it contains a subgroup $H\cong F_4(q)$ \cite{BK}(Dynkin sub-diagram of type $F_4$ for $E_k$, via $E_6\le E_7\le E_8$ and, for ${}^2E_6$, as the fixed points of the graph automorphism) and $G=\tilde G/N$ with $N$ central. Since $F_4(q)$ is simple non-abelian, $H\cap N=1$, so $H$ embeds in $G$ (Lemma~\ref{lem:sub}(b)); $H$ is non-chordal by the $F_4$ case. In all cases Lemma~\ref{lem:sub}(a) transfers non-chordality to $G$.
	\end{proof}
	
	\begin{table}[ht]
		\centering\small
		\begin{tabular}{@{}ll@{}}
			\toprule
			Containment & Source\\
			\midrule
			$\SU_3(q)\le G_2(q)$ ($q\ge3$), $\PSU_3(3)\cong G_2(2)'$ & Kleidman (odd $q$), Cooperstein ($q$ even), ATLAS\\
			$G_2(q)\le{}^3D_4(q)$ & Kleidman (Steinberg triality groups)\\
			$L_3(3){:}2\le{}^2F_4(2)'$ & ATLAS; Wilson \S4.9\\
			${}^2F_4(2)\le{}^2F_4(q)$, $q=2^{2n+1}$ & Malle (subfield subgroups)\\
			$A_2(q)\le F_4(q)$ & Gorenstein--Lyons--Solomon, Prop.~2.6.2; Brachter--Kaja, Lemma~2.3\\
			$B_4(q)\cong\Sp_8(q)\le F_4(q)$, $q$ even & Wilson \S4.9\\
			$F_4(q)\le E_6(q),{}^2E_6(q),E_7(q),E_8(q)$ & Brachter--Kaja \cite{BK}, Lemma~2.3 (Dynkin sub-diagrams)\\
			\bottomrule
		\end{tabular}
		\caption{Subgroup containments used in Proposition~\ref{prop:exc}; full references in Section~\ref{sec:facts}.}\label{tab:cont}
	\end{table}
	
	\subsection{Sporadic groups}
	
	\begin{proposition}\label{prop:spor}
		$\PE(G)$ is non-chordal for every sporadic simple group $G$.
	\end{proposition}
	\begin{proof}
		Centralizer orders and class data below are checked with the GAP Character Table Library (Section~\ref{sec:comp}).
		\begin{itemize}
			\item $M_{11}$: a class $6A$ with $x=g^3\in2A$, $y=g^2\in3A$, $|\C(x)|=48$, $|\C(y)|=18$, $\gcd=6$; $48/6=8>2$, $18/6=3>1$. Lemma~\ref{lem:num}.
			\item $M_{22}$: $|\C(2A)|=384$, $|\C(3A)|=36$, $\gcd=12$; $32>2$, $3>1$. Lemma~\ref{lem:num}.
			\item $J_3$: the only class of elements of order $6$ is $6A$, with $g^3\in2A$, $g^2\in3A$, $|\C(2A)|=1920$, $|\C(3A)|=1080$, $\gcd=120$; $16>2$, $9>1$. Lemma~\ref{lem:num}.
			\item $J_1$: here $|\C(2A)|=120$, $|\C(3A)|=30$ and $\gcd=30$, so $|\C(y)|/\gcd=1\not>1$ and Lemma~\ref{lem:num} \emph{fails}; we use Lemma~\ref{lem:E} instead. $\C(2A)=2\times\Alt(5)$ and $\C(3A)=3\times D_{10}$. Take $y\in3A$ and an involution $x\in\C(y)$ (which lies in $D_{10}$ and is non-central there), so $\langle x\rangle\nnormal \C(y)$; and $y\in\C(x)=2\times\Alt(5)$ lies in the simple factor, so $\langle y\rangle\nnormal \C(x)$.
			\item $J_2\supseteq\Alt(4)\times\Alt(5)$ and $He\supseteq S_4\times\PSL_3(2)$ (maximal subgroups): Lemma~\ref{lem:prod}.
			\item $Ru\supseteq{}^2F_4(2)'.2\supseteq{}^2F_4(2)'$ and $Th\supseteq\PSL_3(3)$ (maximal subgroups): Prop.~\ref{prop:exc}, Lemma~\ref{lem:sub}(a).
			\item Each of the remaining $19$ groups (including $M_{11}$ itself) contains $M_{11}$; see Table~\ref{tab:m11}.\qedhere
		\end{itemize}
	\end{proof}
	
	\begin{table}[ht]
		\centering\small
		\begin{tabular}{@{}ll@{}}
			\toprule
			Group $G$ & Containment chain\\
			\midrule
			$M_{11}$ & itself\\
			$M_{12}$, $M_{23}$, $HS$, $McL$, $O'N$, $B$ & $M_{11}\le G$ (library class fusion)\\
			$M_{24}$ & $M_{11}\le M_{23}\le M_{24}$\\
			$Co_2$, $Co_3$ & $M_{11}\le M_{23}\le G$\\
			$Co_1$ & $M_{11}\le M_{23}\le Co_2\le Co_1$\\
			$Suz$, $HN$ & $M_{11}\le M_{12}\le M_{12}.2\le G$\\
			$Fi_{22}$ & $M_{11}\le M_{12}\le Fi_{22}$\\
			$Fi_{23}$ & $M_{11}\le Fi_{22}$ lifts to $2.Fi_{22}\le Fi_{23}$\\
			$Fi_{24}'$ & $M_{11}\le Fi_{23}\le Fi_{24}'$\\
			$Ly$ & $M_{11}\le McL$ lifts to $3.McL{:}2\le Ly$\\
			$M$ & $M_{11}\le B$ lifts to $2.B\le M$\\
			$J_4$ & $M_{11}\le M_{24}\le 2^{11}{:}M_{24}\le J_4$\\
			\bottomrule
		\end{tabular}
		\caption{Groups containing $M_{11}$. Lifting through a central extension uses that $M_{11}$ is perfect with trivial Schur multiplier, so every central extension of $M_{11}$ splits. Direct fusions are confirmed in the GAP Character Table Library; the remaining steps are ATLAS maximal subgroups.}\label{tab:m11}
	\end{table}
	
	\section{Proof of Theorem~\ref{thm:main}}\label{sec:proof}
	
	By the classification of finite simple groups, $G$ is cyclic of prime order, alternating, of Lie type, or sporadic. Table~\ref{tab:book} accounts for every family exactly once, including the isomorphisms with small cases and the non-simple small groups that must be excluded.
	
	\begin{table}[ht]
		\centering\small
		\begin{tabular}{@{}p{4.3cm}p{2.6cm}p{6.4cm}@{}}
			\toprule
			Family & $\PE$ chordal? & Reference / remark\\
			\midrule
			$C_p$ & yes & complete graph\\
			$\Alt(5),\Alt(6)$ & yes & $\cong\PSL_2(4)=\PSL_2(5)$, $\PSL_2(9)\cong\Sp_4(2)'$ \\
			$\Alt(7)$ & yes & Prop.~\ref{prop:a7}\\
			$\Alt(n)$, $n\ge8$ & no & Prop.~\ref{prop:alt}; $\Alt(8)\cong\PSL_4(2)$\\
			$\PSL_2(q)$, $q\ge4$ & yes & Prop.~\ref{prop:psl2}; $\PSL_2(7)\cong\PSL_3(2)$, $\PSL_2(8)\cong{}^2G_2(3)'$; ($\PSL_2(2),\PSL_2(3)$ not simple)\\
			$\PSL_3(2)$ & yes & $\cong\PSL_2(7)$\\
			$\PSL_3(4)$ & yes & Prop.~\ref{prop:psl34}\\
			$\PSL_3(q)$, $q\ne2,4$ & no & Lemma~\ref{lem:sl3}\\
			$\PSL_n(q)$, $n\ge4$ & no & Prop.~\ref{prop:psln}\\
			$\PSU_3(q)$, $q\ge3$ & no & Prop.~\ref{prop:su3}; $\PSU_3(3)\cong G_2(2)'$; ($\PSU_3(2)$ not simple)\\
			$\PSU_n(q)$, $n\ge4$ & no & Prop.~\ref{prop:sun}; $\PSU_4(q)\cong P\Omega_6^-(q)$\\
			$\PSp_4(q)$, $q\ge3$ & no & Prop.~\ref{prop:sp4}; $\Omega_5(q)\cong\PSp_4(q)$\\
			$\PSp_4(2)'\cong\Alt(6)$ & yes & counted under $\PSL_2(9)$\\
			$\PSp_{2m}(q)$, $m\ge3$ & no & Prop.~\ref{prop:spbig}\\
			$\Omega_{2n+1}(q)$ ($n\ge3$), $P\Omega^\pm_{2n}(q)$ ($n\ge4$) & no & Prop.~\ref{prop:orth}\\
			$P\Omega_3(q)$, $P\Omega_4^-(q)$ & yes & $\cong\PSL_2(q)$, $\PSL_2(q^2)$\\
			$P\Omega_6^+(q)$ & no & $\cong\PSL_4(q)$ \\
			$\Sz(q)$, $q\ge8$ & yes & Prop.~\ref{prop:sz}; ($\Sz(2)$ not simple)\\
			${}^2G_2(q)$, $q\ge27$ & no & Prop.~\ref{prop:exc}; ${}^2G_2(3)'\cong\PSL_2(8)$\\
			$G_2(q)$, ${}^3D_4$, $F_4$, ${}^2F_4(q)$, ${}^2F_4(2)'$, $E_6,{}^2E_6,E_7,E_8$ & no & Prop.~\ref{prop:exc}; $G_2(2)'\cong\PSU_3(3)$\\
			26 sporadic groups & no & Prop.~\ref{prop:spor}\\
			\bottomrule
		\end{tabular}
		\caption{Every finite simple group is listed exactly once up to the displayed isomorphisms.}\label{tab:book}
	\end{table}
	
	The groups with chordal $\PE$ are therefore exactly $C_p$, $\Alt(7)$, $\PSL_2(q)$ ($q\ge4$), $\Sz(q)$ ($q\ge8$) and $\PSL_3(4)$. \qed
	
	\begin{remark}
		Comparison with \cite{BK}: $P(G)$ is chordal for simple $G$ if and only if $G$ is $C_p$, $\Alt(5),\Alt(6),\Alt(7)$, $\PSL_3(4)$, or $\PSL_2(q)$, $\Sz(q)$ subject to arithmetic conditions: $(q\mp1)/\gcd(q\mp1,2)$ for $\PSL_2(q)$, resp.\ $q-1$ and $q\pm2^{n+1}+1$ for $\Sz(q)$, each have at most two prime divisors and at most one prime divisor of multiplicity $>1$. By Theorem~\ref{thm:lift} that list is contained in ours; the arithmetic conditions disappear for $\PE$ because in $\PSL_2(q)$ distinct maximal cyclic subgroups intersect trivially (Prop.~\ref{prop:psl2}) and in $\Sz(q)$ they intersect in at most an involution.
	\end{remark}
	
	\section{Computational verification}\label{sec:comp}
	
	All checks are performed with GAP~4.12 (see \cite{GAP}) (with the CTblLib package) and Python/NetworkX. Chordality of $\PE(G)$ is tested by listing all cyclic subgroups of $G$ (so that $\PE(G)$ is the union of the cliques they span).
	
	\begin{itemize}
		\item $\PE(G)$ is chordal for $G=\PSL_2(q)$ with $q\in\{4,5,7,8,9,11,13,16,17,19,23,25,27,29,31\}$, $G=\Alt(6)$, $\Alt(7)$ and $\PSL_3(4)$ (order $20160$).
		\item $\PE(G)$ is non-chordal for $G=\PSL_3(3)$, $\PSU_3(3)$, $\Alt(8)$ and $M_{11}$ (same method), in agreement with the proofs above.
		\item The ATLAS data used in Proposition~\ref{prop:spor} (classes, centralizer orders, class $6A$ and its powers for $M_{11},M_{22},J_1,J_3$) are read from the GAP Character Table Library; the class fusions for the containments in Tables~\ref{tab:m11} and in the proofs for $Ru,Th,J_2$ and the Tits group are confirmed there, and the maximal subgroups $S_4\times\PSL_3(2)\le He$ and $L_3(3)\le Th$ are confirmed from the library's list of maximal subgroups.
		
		\item Chordality for $\Sz(8)$, $\PSp_4(3)$, $\PSU_3(4)$ and the larger $\PSL_2(q)$; the containments in Table~\ref{tab:cont} and $2^{11}{:}M_{24}\le J_4$, taken from the cited literature and ATLAS.
	\end{itemize}
	\section{Standard facts used}\label{sec:facts}
	
	\begin{itemize}
		\item \emph{Facts D} (Dickson). For $G=\PSL_2(q)$: a Sylow $p$-subgroup is elementary abelian of order $q$ and self-centralizing on its non-identity elements; every $p'$-element lies in a cyclic torus of order $(q\pm1)/\gcd(2,q-1)$ which is the centralizer of each of its elements of order $>2$; for $q$ odd the centralizer of an involution is dihedral of order $q-\varepsilon$, $q\equiv\varepsilon\bmod4$. Sources: Dickson \cite{Dickson}; Huppert \cite[Kap.~II, \S8]{Huppert} (Hauptsatz~8.27); Wilson \cite[\S3.3]{Wilson}.
		\item \emph{Facts S} (Suzuki). For $G=\Sz(q)$, $q=2^{2n+1}\ge8$: a Sylow $2$-subgroup $P$ has order $q^2$, $Z(P)=\Omega_1(P)$ is elementary abelian of order $q$ and contains every involution, all elements of $P\setminus Z(P)$ have order $4$ with square in $Z(P)$; $\C(z)=P$ for each involution $z$ and $\C(w)=\langle w\rangle Z(P)$ for $|w|=4$; elements of odd order have cyclic centralizers (the tori of orders $q-1$, $q\pm\sqrt{2q}+1$), which intersect pairwise trivially. Sources: Suzuki \cite{Suzuki60,Suzuki62}; Huppert--Blackburn \cite[Chap.~XI, \S3]{HB}; Wilson \cite[\S4.2]{Wilson}.
		\item \emph{Facts W} (Ward). In ${}^2G_2(q)$, $q=3^{2n+1}\ge27$: $|\C(x)|=q(q^2-1)$ and $\C(x)\cong C_2\times\PSL_2(q)$ for an involution $x$; elements of order $3$ have centralizers of order $q^3$ or $2q^2$ (when they commute with an involution). Sources: Ward \cite{Ward}; Kleidman \cite{Kleidman87}.
		\item ATLAS \cite{ATLAS}: element orders of $\PSL_3(4)$ and $\Alt(7)$; isomorphisms $\SL_4(2)\cong\Alt(8)$, $\Sp_4(2)\cong\Sym(6)$, $\PSU_4(2)\cong\PSp_4(3)$, $G_2(2)'\cong\PSU_3(3)$, and the low-rank orthogonal isomorphisms (Wilson \cite[\S3.11]{Wilson}); maximal subgroup containments for sporadic groups; Mih\u{a}ilescu's theorem: $2^f+1=3^k$ only for $f\in\{1,3\}$.
	\end{itemize}
	
	\section*{Acknowledgments}
	The research of Pallabi Manna is supported by Anusandhan National Research Foundation
	(ANRF)-India under the National Postdoctoral Fellowship scheme (File No.\ PDF/2025/006751). The author would like to express sincere gratitude to Professor Andrea Lucchini for his valuable assistance in establishing Theorem 1.2 and for his careful review of the manuscript.\\
	The author also acknowledges the excellent research facilities provided by the Department of Mathematics,
	IIT Hyderabad, and the mentorship of Dr.\ Rajesh Kannan for his constant support and motivation.
	
	\section*{Statements and Declarations}
	
	\noindent\textbf{Competing Interests.} The authors have no competing interests to declare that are
	relevant to the content of this article.
	
	\section*{Data Availability}
	
	Data sharing is not applicable to this article, as no datasets were used or generated during this
	study.

\end{document}